\documentclass[a4paper,12pt]{amsart}

\usepackage{amsmath, amssymb}

\newtheorem{lemma}{Lemma}[section]
\newtheorem{corollary}[lemma]{Corollary}
\newtheorem{theorem}[lemma]{Theorem}

\newtheorem{remark}{Remark}
\newcounter{claim}[lemma]

\renewcommand{\unlhd}{\,\underline{\triangleleft}\,}
\def \udot  {{}^{\textstyle .}}
\def \Syl {\mathrm {Syl}}
\def \syl {\mathrm {Syl}}
\def \Sym {\mathrm \Sym}
\def \Dih {Dih}

\def \cal {\mathcal}

\def \AA {{\cal A}}

\def \e {\eta}

\def \G {\gamma}

\def \Inn{\mbox {\rm Inn}}

\def \Stab {\mbox {\rm Stab}}

\def \Sym {\mbox {\rm Sym}}

\def \Im {\mbox {\rm Im}}
\def \PSU {\mbox {\rm PSU}}
\def \SU {\mbox {\rm SU}}
\def \GU {\mbox {\rm GU}}

\def \PGammaL {\mbox {\rm P}\Gamma {\rm L}}
\def \PSigmaL {\mbox {\rm P}\Sigma {\rm L}}
\def \PGL {\mbox {\rm PGL}}

\def \Dih {\mbox {\rm Dih}}

\def \SL {\hbox {\rm SL}}

\def \GL {\mbox {\rm GL}}
\def \syl {\hbox {\rm Syl}}

\def \Aut{ \mathrm {Aut}}

\def \J{\mbox {\rm J}}

\def \M{\mbox {\rm M}}

\def \udot {{}^{\textstyle .}}

\def \DD   {\mbox {{\rm D}}}

\def \U {\mbox {{\rm U}}}

\def \Alt {\mbox {{\rm Alt}}}
\def \GG {\mbox {{\rm G}}}

\def \K {{\mathrm {K}}}

\def \Q {\mbox {{\rm Q}}}

\def \PSL {\mbox {{\rm PSL}}}

\def \normal {\unlhd}
\def \G {\Gamma}

\def \e {\eta}

\begin{document}

\title{Semisymmetric cubic graphs of twice odd order}
\author{Chris Parker}
\date{\today}
\address{Chris Parker\\
School of Mathematics and Statistics\\
University of Birmingham\\ Edgbaston\\ Birmingham B15 2TT, UK}
\email{C.W.Parker@bham.ac.uk}
\thanks{{\em Mathematical subject classification}: 20D08}
\maketitle

\section{Introduction}

Suppose that $\Gamma$ is a connected graph and $G$ is a subgroup
of the automorphism group $\Aut(\Gamma)$ of  $\Gamma$. Then
$\Gamma$ is $G$-\emph{symmetric} if $G$ acts transitively on the
arcs (and so the vertices) of $\Gamma$ and $\Gamma$ is
$G$-\emph{semisymmetric} if $G$ acts edge transitively but
\emph{not} vertex transitively on $\Gamma$. If $\Gamma$ is
$\Aut(\Gamma)$-symmetric, respectively,
$\Aut(\Gamma)$-semisymmetric, then we say that $\Gamma$ is
\emph{symmetric}, respectively, \emph{semisymmetric}.  If $\Gamma$
is $G$-semisymmetric, then we say that $G$ acts
\emph{semisymmetrically} on $\Gamma$. If $\Gamma$ is
$G$-semisymmetric, then the orbits of $G$ on the vertices of
$\Gamma$ are the two parts of a bipartition of  $\Gamma$.

Semisymmetric cubic graphs (graphs in which every vertex has
degree 3) have been the focus of a number of recent  articles, we
mention specifically
\cite{Du1,Du12,Du3,Lip,Lu,Malnic,Malnic2,Marusic} where infinite
families of such graphs are presented and where the semisymmetric
graphs of order $2pq$, $2p^3$, $6p^2$ with $p$ and $q$ odd primes
 are determined (with the help of the classification
of the finite simple groups). We also remark that a catalogue of
all the semisymmetric cubic graphs of order at most 768 has
recently been obtained by Conder \emph{et. al.} \cite {Conder}.
The objective of this article is to partially describe all groups
which act semisymmetrically on a   cubic graph of order twice an
odd number. Thus our result reduces some of the aforementioned
investigations to checking group orders in our list. We also
mention that our theorems call upon only a small number of
characterization theorems used in the classification of the finite
simple groups and not on the whole classification itself.

Suppose that $\Gamma$ is a $G$-semisymmetric cubic graph. Let
$\{u,v\}$ be an edge in $\Gamma$. Set $G_u = \Stab_G(u)$, $G_v =
\Stab_G(v)$ and $G_{uv}= G_u \cap G_v$. Then, as $G$ acts edge
transitively on $\G$ and $u$ is not in the same $G$-orbit as $v$,
we have $[G_u:G_{uv}]= [G_v:G_{uv}]=3$. Suppose that $K \normal G$
and $K \le G_{uv}$. Then $K$ fixes every edge of $\Gamma$ and
hence $K=1$. As $\Gamma$ is connected, the subgroup $\langle
G_u,G_v\rangle$ acts transitively on the edges of $\Gamma$ and so
we infer that $G=\langle G_u,G_v\rangle$. We have shown that $G$
satisfies
\begin{itemize}
\item  $G = \langle G_u,G_v\rangle$;
 \item
$[G_u:G_u\cap G_v]=[G_v:G_u\cap G_v]=3$; and \item no non-trivial
subgroup of $G_{uv}$ is normal in $G$.
\end{itemize}
This group theoretic configuration has been studied by Goldschmidt
in \cite{Goldschmidt} where it is shown that  the triple
$(G_u,G_v,G_{uv})$ is isomorphic (as an amalgam) to one of fifteen
possible such triples (see Table~\ref{GoldList}). Thus if $\Gamma$
is $G$-semisymmetric, then the structures of $G_u$, $G_v$ and
$G_{uv}$ (and the embeddings of $G_{uv}$ into $G_u$ and $G_v$) are
known (up to swapping the roles of $u$ and $v$). We call the
possible triples of groups appearing in Table~\ref{GoldList}
\emph{Goldschmidt amalgams}. To bring more precision to the
definition of an  amalgam, we define an \emph{amalgam} to be a
pair of monomorphisms $(\phi_u:G_{uv}\rightarrow G_u,
\phi_v:G_{uv}\rightarrow G_v)$ and then define a group $G$ to be a
\emph{completion} of $(\phi_u,\phi_v)$ provided there exist
monomorphisms $\psi_u: G_u \rightarrow G, \psi_v:G_v\rightarrow G$
such that $\phi_u\psi_u=\phi_v\psi_v$ and $G =\langle
\Im(\psi_u),\Im(\psi_v)\rangle$. In the case that $G$ is a
completion of an amalgam, we identify the groups $G_u$, $G_v$ and
$G_{uv}$ with their images in $G$. Furthermore, when there is no
danger of confusion, we identify the amalgam $(\phi_u,\phi_v)$
with the triple of subgroups  $(G_u,G_v,G_{uv})$ and the
monomorphisms $\phi_u$ and $\phi_v$ are then understood to be
inclusion maps. In the case of the Goldschmidt amalgams, this
leads to no confusion as the amalgams are determined uniquely by
the structure of the pairs of subgroups $G_u$ and $G_v$ together
with the fact that $G_{uv}$ contains no non-trivial subgroup $K$
with the property that $\phi_u(K)\normal G_u$ and $\phi_v(K)
\normal G_v$. Assume that $G$ is a completion of an amalgam $\AA=
(G_u,G_v,G_{uv})$. We construct  a graph $\G =\Gamma(G,\AA)$ which
has vertex set $\{G_ug, G_v g \mid g \in G\}$ and edges
$\{\{G_uh,G_vg\}\mid G_uh\cap G_vg\not=\emptyset\}$. The graph
$\Gamma$ is called the \emph{coset graph} of the completion $G$.
The group $G$ acts on $\G$ by right multiplication. The kernel $K$
of this action is the largest normal subgroup of $G$ which is
contained in $G_{uv}$. The graph $\Gamma$  is then
$G/K$-semisymmetric. In case $\AA$ is a Goldschmidt amalgam, we
have that $K=1$ and so $\Gamma$ is $G$-semisymmetric. So
understanding $G$-semisymmetric cubic graphs is the same as
understanding completions of  Goldschmidt amalgams. We mention
that in \cite{ParkerRowley1,ParkerRowley2} completions of the
Goldschmidt amalgam of type $\GG_3$ in $\PSL_3(p^a)$  and
completions which are sporadic simple groups are investigated (see
Table~\ref{GoldList} for an explanation of the notation). In
\cite{ParkerRowley3} completions of the Goldschmidt amalgam of
type  $\GG_4$ which are contained in $\PSL_3(p^a)$ are determined.
We shall call upon this result later in this paper. For the
investigations in this article, we are interested in finite groups
$G$ which are completions of a Goldschmidt amalgam $\AA$ and which
have $[G:G_u]+[G:G_v] = 2[G:G_v]$ equal to twice an odd number.
Since, by Goldschmidt's Theorem, $G_{uv}$ is a $2$-group, this
means that $G_{uv}$ is a Sylow $2$-subgroup of $G$. If $G$ is  a
completion of a Goldschmidt amalgam and $G_{uv} \in \Syl_2(G)$,
then we call $G$
 a \emph{Sylow completion} of the amalgam $\AA$.
We emphasise that, so long as $G \not=G_1$,  the graph
$\Gamma(G,\AA)$ corresponding to a Sylow completion $G$ of $\AA$
is a  $G$-semisymmetric cubic graph of twice odd order. We need
one final definition before we can state our main theorem. Suppose
that $G$ is a  completion of an amalgam $\AA$. Then a normal
subgroup $R$ of $G$ is called a \emph{regular} normal subgroup of
$G$ provided $R$ acts semiregularly on the vertices of
$\Gamma(G,\AA)$.

\begin{theorem}\label{MainTheorem} Suppose that $G$ is a Sylow
completion of a Goldschmidt amalgam.
 Then there exists
a regular normal subgroup $R$ of $G$ of odd order such that $G/R$
is isomorphic to one of the groups  in column three of
Table~\ref{Answers}.
\end{theorem}

We also give the graph theoretic version of
Theorem~\ref{MainTheorem}.

\begin{theorem}\label{MainTheorem2} Suppose that $G$ acts semisymmetrically on a cubic
graph $\Gamma$  of twice odd order.
 Then there exists
a normal subgroup $R$ of $G$ of odd order which acts semiregularly
on the vertices of $\Gamma$  and such that  $G/R$ is isomorphic to
one of the groups  in column three of Table~\ref{Answers}.
\end{theorem}

\begin{table}
\begin {center}
\begin{tabular}{||cccl||}
\hline
 \hline
Division& Type&$G/R$&Condition\\
 \hline
1& $\GG_1$& $3$&\\
&$\GG_1^1$&$\Sym(3)$&\\
\hline

2&$\GG_1$ &$3^2$&\\
&$\GG_1^1$&$3^2.2$&\\
&$\GG_1^2$&$3\wr 2$&\\
&$\GG_1^3$&$\Sym(3)\times \Sym(3)$&\\

\hline 3&
$\GG_2$&$3^3.\Alt(4)$&\\
&$\GG_2^1$&$3^3.\Sym(4)$&\\
&$\GG_2^2$&$3^3.\Sym(4)$&\\
&$\GG_2^3$&$\Sym(3)\wr 3$&\\
&$\GG_2^4$&$\Sym(3)\wr \Sym(3)$&\\
\hline
4&$\GG_1^3$&$\PSL_2(p)$&$p$ a prime, $p\equiv 11,13 \pmod{24}$\\
\hline 5& $\GG_2$&$\PSL_2(p)$&$p$ a prime, $p\equiv 11,13 \pmod{24}$\\
& $\GG_2^1$&$\PGL_2(p)$&$p$ a prime, $p\equiv 11,13 \pmod{24}$\\
\hline 6& $\GG_2^2$&$\Alt(7)$&\\
&$\GG_2^4$&$\Sym(7)$&\\
\hline 7&$\GG_2^1$&$\PSL_2(p)$&$p$ a prime, $p\equiv 23,25\pmod
{48}$\\\hline 8&
$\GG_2^1$&$\PSL_2(p^2)$&$p$ a prime, $p\equiv 5,19\pmod {24}$\\
&$\GG_2^4$&$\PSigmaL_2(p^2)$&$p$ a prime, $p\equiv 5,19\pmod {24}$\\
\hline 9& $\GG_3$& $\PSL_2(p)$&$p$ a prime, $p\equiv 7,9\pmod
{16}$\\\hline
10& $\GG_3$& $\PSL_2(p^2)$&$p$ a prime, $p\equiv 3,5\pmod {8}$\\
 & $\GG_3^1$& $\PSigmaL_2(p^2)$&$p$ a prime, $p\equiv 3,5\pmod {8}$\\
  \hline
11& $\GG_4$& $\PSL_3(p)$&$p$ a prime, $p\equiv 5\pmod {8}$\\
  &$\GG_4^1$& $\PSL_3(p).2$&$p$ a prime, $p\equiv 5\pmod {8}$\\
\hline
12& $\GG_4$& $\PSU_3(p)$&$p$ a prime, $p\equiv 3\pmod {8}$\\
  &$\GG_4^1$& $\PSU_3(p).2$&$p$ a prime, $p\equiv 3\pmod {8}$\\
\hline
  13&$\GG_5$& $\M_{12}$&\\
  &$\GG_5^1$& $\Aut(\M_{12})$&\\
  \hline
14&$\GG_5$& $\GG_2(p)$&$p$ a prime, $p\equiv 3,5\pmod {8}$\\
&$\GG_5^1$& $\Aut(\GG_2(3))$&\\
\hline \hline
\end{tabular}
\end{center}
\caption{Sylow completions of Goldschmidt amalgams}\label{Answers}
\end{table}

Suppose that $G$ is a completion of a Goldschmidt amalgam
$\AA=(G_u,G_v,G_{uv})$, $\G =\G(G,\AA)$ and $R$ is a regular
normal subgroup of $G$. Define $\Gamma/R$ to be the coset graph of
the amalgam $(G_vR, G_uR, G_{uv}R)$. Then $\Gamma/R$ is admits $G$
by right multiplication and  $R$ in the kernel. Furthermore,
unless $G_uR = G_vR =G$, $\Gamma/R$ is a $G/R$-semisymmetric cubic
graph. The quotients $G/R$ in division 1 of Table~\ref{Answers} do
not themselves determine interesting graphs as in this case
$\Gamma/R$ has only two vertices. On the other hand,  division 1
may by no means be omitted from the table as is shown by the
semisymmetric graphs  {\bf S486} and {\bf S702c} of \cite{Conder}.

The groups $G/R$ in division $2$ of Table~\ref{Answers} all have
$\Gamma/R$ isomorphic to the symmetric graph $\K_{3,3}$. However,
the group $7^2.(3\wr 2)$ is a Sylow completion of an amalgam of
type $\GG_1^2$ and gives the semisymmetric graph listed as {\bf
S294} in \cite{Conder} and so we see that there are also
semisymmetric graphs appearing in this division. All the groups
$G/R$ in division $3$ are acting on the smallest semisymmetric
graph, the Gray graph, listed as {\bf S54} in \cite{Conder} (see
also \cite{Bouwer}).

The graphs arising from the groups in division 7 of
Table~\ref{Answers} are related to the those investigated by
Lipschutz and Xu  in \cite{Lip}. The author expects that all the
 cubic $G$-semisymmetric graphs arising from
$\PSL_2(p^a)$ are comparatively well-known. On the other hand, the
graphs constructed from the groups in divisions 11, 12 and 14 of
Table~\ref{Answers} appear to be new. Finally we remark that the
coset graph constructed from $\Aut(\M_{12})$ (division 13) was
investigated a number of years ago by Biggs \cite{Biggs}.

Of course, even though $G/R$  acts semisymmetrically on
$\Gamma/R$, it may be that $\Gamma/R$ is nonetheless a symmetric
graph as $\Aut(\Gamma/R)$ may act transitively on the vertices of
$\Gamma/R$. For $G$ and $R$ as in Table~\ref{Answers}, we
investigate the automorphism groups of $\Gamma/R$ and whether or
not the graph is symmetric or semisymmetric. The result is
reported in the following corollary.

\begin{corollary}\label{gammaR} Suppose that $\Gamma$, $G$ and $R$
are as in Theorem~\ref{MainTheorem2}. Then, for the cases in
divisions 2 to 14 in Table~\ref{Answers},  $\Aut(\Gamma/R)$ and
whether or not  $\G/R$ is symmetric or semisymmetric is listed in
Table~\ref{SS}.
\end{corollary}

\begin{table}
\begin {center}
\begin{tabular}{||cccc||}
\hline
 \hline
Division& $\Aut(\G/R)$&Condition &Symmetric\\
 \hline
2&$\Sym(3)\wr 2$ &&Yes\\\hline
 3&$\Sym(3)\wr \Sym(3)$&&No \\
 \hline
4&$\PGL_2(p)$&$p$ a prime, $p\equiv 11,13 \pmod{24}$&Yes\\
\hline 5&$\PGL_2(p)$&$p$ a prime, $p\equiv 11,13 \pmod{24}$&No\\
\hline 6&$\Sym(7)$&&No\\
\hline 7&$\PSL_2(p)$&$p$ a prime, $p\equiv 23,25\pmod
{48}$&No\\\hline 8&$\PSigmaL_2(p^2)$&$p$ a prime, $p\equiv 5,19\pmod {24}$&No\\
\hline 9&$\PGL_2(p)$&$p$ a prime, $p\equiv 7,9\pmod
{16}$&Yes\\\hline 10&
  $\PGammaL_2(p^2)$&$p$ a prime, $p\equiv 3,5\pmod {8}$&Yes\\
  \hline
11&$\PSL_3(p).2$&$p$ a prime, $p\equiv 5\pmod {8}$&No\\
\hline 12
  & $\PSU_3(p).2$&$p$ a prime, $p\equiv 3\pmod {8}$&No\\
\hline
  13
  & $\Aut(\M_{12})$&&No\\
  \hline
14& $\GG_2(p)$&$p$ a prime, $p\equiv 3,5\pmod {8}$, $p\ge 5$&No\\
&$\Aut(\GG_2(3))$&&No\\
\hline \hline
\end{tabular}
\end{center}
\caption{The automorphism group of $\G/R$}\label{SS}
\end{table}

We can use Theorem~\ref{MainTheorem} to recover the following
result (slightly corrected) of Iofinova and Ivanov
 \cite{IvanovIofinova}.

\begin{corollary}\label{biprimitive} Suppose that $\Gamma$
is a  semisymmetric, cubic  graph
 and that $\Aut(\G)$ acts primitively on both parts of the
bipartition of $\G$. Then $\Aut(\Gamma)$ is isomorphic to one of
$\PGL_2(11)$, $\PGL_2(13)$, $\PSL_2(23)$, $\GG_2(2)$ or
$\Aut(\M_{12})$.
\end{corollary}

In a similar spirt, and at the same time,  we can also determine
those  cubic semisymmetric graphs $\Gamma$ of twice odd order on
which $\Aut(\Gamma)$ acts primitively on exactly one part of the
bipartition of $\Gamma$.

\begin{corollary}\label{uniprimitive} Suppose that $\Gamma$ is a  semisymmetric, cubic  graph of twice odd order and that $\Aut(\G)$ acts primitively on exactly one part of the
bipartition of $\G$. Then one of the following holds.
\begin{enumerate}
\item (division 3, Table~\ref{Answers}) $\Aut(\Gamma) \cong \Sym(3)\wr
\Sym(3)$.
\item (division 5, Table~\ref{Answers}) $\Aut(\Gamma) \cong \PGL_2(p)$, $p\equiv 11,13
\pmod{24}$ and $p > 13$.
\item (division 7, Table~\ref{Answers}) $\Aut(\Gamma) \cong \PGL_2(p)$, $p\equiv
23,15 \pmod{48}$ and $p > 23$.
\item (division 8, Table~\ref{Answers}) $\Aut(\Gamma) \cong \PSigmaL_2(25)$.
\item (division 8, Table~\ref{Answers}) $\Aut(\Gamma) \cong \PSigmaL_2(p)$, $p>5$, $p\equiv 5,19 \pmod {24}$.
\item (division 11, Table~\ref{Answers}) $\Aut(\Gamma) \cong \Aut(\PSL_3(5))$.
\end{enumerate}
\end{corollary}

The paper is organized as follows. In Section~\ref{Sec2} we
present Goldschmidt's fundamental  theorem  and  gather some
results about the Goldschmidt amalgams. In Section 3 we quote some
results which characterise various simple groups by their Sylow
$2$-structure  and in Section 4 we use the results of the previous
two sections to prove Theorem~\ref{MainTheorem} and its
corollaries.

Our group theoretic notation will for the most part follow
\cite{Gorenstein}. We use {\sc Atlas} \cite{Atlas} notation for
group extensions and for describing the shape of groups. The
symmetric group on $n$ letters is denoted by $\Sym(n)$, the
alternating group on $n$ letters is denoted by $\Alt(n)$, the
dihedral group of order $n$ is denoted by $\Dih(n)$ and the
quaternion group of order $8$ is denoted by $\Q_8$. We shall use
both $\mathbb Z_n$ and $n$ to denote the cyclic group of order
$n$.  The remaining notation for the classical and Lie type groups
is standard and hopefully self explanatory. Finally, if a group
$G$ acts on a $p$-group $Q$, then $\e(G,Q)$ denotes the number of
non-central $G$-chief factors in $Q$.

\section{The Goldschmidt amalgams}\label{Sec2}

First we state the classification theorem due to Goldschmidt
\cite{Goldschmidt} which, as mentioned in the introduction, is our
way into the study of semisymmetric graphs.

\begin{theorem}\label{Goldschmidt} Suppose that $(G_1,G_2,G_{12})$ is
an amalgam which satisfies
\begin{enumerate}
\item  $G_i$ is finite and $[G_i:G_{12} ]=3$, for
$i=1,2$; and \item no non-trivial subgroup of $G_{12}$ is normal
in both $G_1$ and $G_2$.
\end{enumerate} Then $|G_{12}|$ divides  $2^7$, $\e(G_i,O_2(G_i)) \le
2$ for $i=1,2$ and up to isomorphism $(G_1,G_2)$ is isomorphic to
one of the pairs of groups listed in
Table~\ref{GoldList}.\end{theorem}

\begin{table}
\noindent \begin{tabular}{||ccl||} \hline\hline
Name&$(G_1,G_2)$&Completion\\
\hline  $\GG_1$&$({\mathbb Z}_3,{\mathbb Z}_3)$&
${\mathbb Z}_3\times {\mathbb Z}_3$\\
$\GG_1^1$&$(\Sym(3),\Sym(3))$&
$3^2.2$\\
$\GG_1^2$&$(\Sym(3),{\mathbb Z}_6)$&
$\Sym(3) \times {\mathbb Z}_3$\\
$\GG_1^3$&$(\Sym(3)\times {\mathbb Z}_2,\Sym(3)\times {\mathbb
Z}_2)$&
$\PSL_2(11)$\\

$\GG_2$&$(\Alt(4),\Sym(3)\times {\mathbb Z}_2)$&$\PSL_2(11)$\\
$\GG_2^1$&$(\Sym(4), \Dih(24))$&$\PSL_2(23)$\\
$\GG_2^2$&$(\Sym(4),(2^2\times 3).2) ) $&$\Alt(7)$\\
$\GG_2^3$&$(\Alt(4) \times {\mathbb Z}_2, \Sym(3) \times {\mathbb Z}_2 \times {\mathbb Z}_2)$&$\Sym(3) \wr {\mathbb Z}_3$\\
$\GG_2^4$&$(\Sym(4) \times {\mathbb Z}_2, \Sym(3)
\times \Dih(8))$&$\Sym(7)$\\

$\GG_3$&$(\Sym(4),\Sym(4))$&$\Alt(6)$, $\SL_3(2)$\\
$\GG_3^1$&$(\Sym(4) \times {\mathbb Z}_2,\Sym(4)\times {\mathbb Z}_2)$&$\Sym(6)$\\
$\GG_4$&$(({\mathbb Z}_4 \times {\mathbb Z}_4).\Sym(3), (\Q_8 \circ {\mathbb Z}_4).\Sym(3))$&$\GG_2(2)^\prime \cong \U_3(3)$\\
$\GG_4^1$&$(({\mathbb Z}_4 \times {\mathbb Z}_4).(\Sym(3)\times {\mathbb Z}_2), (\Q_8 \circ \Q_8).\Sym(3)),$&$\GG_2(2)$\\
& $\;\;\;\;\;\;(\e(G_2,O_2(G_2))=1$)&\\
$\GG_5$&$(({\mathbb Z}_4 \times {\mathbb Z}_4).(\Sym(3)\times {\mathbb Z}_2), (\Q_8 \circ \Q_8).\Sym(3)),$&$\M_{12}$\\
& $\;\;\;\;\;\;(\e(G_2,O_2(G_2))=2$)&\\
 $\GG_5^1$&$(({\mathbb Z}_4 \times {\mathbb Z}_4).(2^2 \times 3).2, (\Q_8 \circ \Q_8).(\Sym(3) \times{\mathbb Z}_2)),$&$\Aut(\M_{12})$\\
& $\;\;\;\;\;\;(\e(G_2,O_2(G_2))=2$)&\\
\hline\hline
\end{tabular}
 \caption{The Goldschmidt amalgams}\label{GoldList}
\end{table}
\medskip

The sample completions presented in Table~\ref{GoldList} tell us
exactly the isomorphism types of the groups $G_1, G_2$ and
$G_{12}$ in a Goldschmidt amalgam. In particular, we emphasize the
following. Let $\mathcal A =(G_1,G_2,G_{12})$ be a Goldschmidt
amalgam. Then, if $\mathcal A$ has type $\GG_4$, $G_{12}$ is
isomorphic to a Sylow $2$-subgroup of $\PSU_3(3)$, if $\mathcal A$
has type $\GG_4^1$ or $\GG_5$, then $G_{12}$ is isomorphic to a
Sylow $2$-subgroup of $\M_{12}$. (To see this we notice that
$\GG_2(\mathbb Z) \cong \GG_2(2) $ and this is contained as a
subgroup of $\GG_2(3)$. Since $\GG_2(3)$ and $\M_{12}$ have
isomorphic Sylow $2$-subgroups the claim follows.) Finally, if
$\mathcal A$ is of type $\GG_5^1$, then $G_{12}$ is isomorphic to
a Sylow $2$-subgroup of $\Aut(\M_{12})$.

The \emph{type} of a Goldschmidt amalgam is simply the name given
to it in column one of Table~\ref{GoldList}. The \emph{class} of a
Goldschmidt amalgam of type $\GG_i^j$ is $\GG_i$. So, for example,
the amalgams $\GG_1$, $\GG_1^1$, $\GG_1^2$, and $\GG_1^3$ are all
in Goldschmidt class $\GG_1$.

Notice also that if $(G_1,G_2,G_{12})$ is a Goldschmidt amalgam
then so is $(G_2,G_1,G_{12})$ and that these amalgams only appear
once in Table~\ref{GoldList}. The reader should be cautious when
reading the structure of an amalgam from the list as we may have
swapped the ordering of $G_1$ and $G_2$. This point is made more
clear in   the next lemma.

\begin{lemma}\label{gold1} Suppose that $\mathcal A=(G_1,G_2,G_{12})$ is a
Goldschmidt amalgam.
\begin{enumerate}
\item If $\mathcal A$ is in class $\GG_1$, then, for $i\in\{1,2\}$,
$[O_2(G_i),O^2(G_i)]=1$.

\item If $\mathcal A$ is in class $\GG_2$, then there exists  $i\in
\{1,2\}$ such that $[O_2(G_i),O^2(G_i)]=1$ and
$[O_2(G_{3-i}),O^2(G_{3-i})]\not=1$.

\item If $\mathcal A$ is in class $\GG_3$, $\GG_4$ or $\GG_5$, then
$[O_2(G_1),O^2(G_1)]\not=1$ and $[O_2(G_2),O^2(G_2)]\not=1$.

 \item If $\mathcal A$ is in class $\GG_3$, $\GG_4$ or $\GG_5$,
then, for $i\in \{1,2\}$, $O_2(O^2(G_i)) \not\le O^2(G_{3-i})$.

\end{enumerate}
\end{lemma}

\begin{proof} All the claims  follow from the structure of the Goldschmidt amalgams
given in Table~\ref{GoldList}.
\end{proof}

The next lemma shows that an amalgam of Goldschmidt class $\GG_k$
contains a subamalgam of Goldschmidt type $\GG_k$.

\begin{lemma}\label{subamalgam} Suppose that $(G_1,G_2,G_{12})$ is a Goldschmidt
amalgam of type $\GG_k^j$. For $i=1,2$, let
$M_i=O^2(G_i)O_2(O^2(G_{3-i}))$ and set $M_{12} = M_1\cap M_2
=O_2(O^2(G_1))O_2(O^2(G_2))$. Then $(M_1,M_2,M_{12})$ is  a
Goldschmidt amalgam of type $\GG_i$.
\end{lemma}

\begin{proof}  We have that  $M_{12}$
is a Sylow $2$-subgroup of both $M_1$ and $M_2$. If $K \le M_{12}$
is chosen of maximal order subject to being  normalized by both
$M_1$ and $M_2$, then, as $G_{12}$ normalizes both $M_{1}$ and
$M_2$, $K^x=K$ for all $x \in G_{12}$. Hence $K=1$ as
$(G_1,G_2,G_{12})$ is a Goldschmidt amalgam. It is now elementary
to check that $(M_1,M_2,M_{12})$ is of type $\GG_k$.
\end{proof}

We now introduce some important notation.  For $G$ a completion of
the amalgam $(G_1,G_2,G_{12})$ and $N$ a normal subgroup of $G$.
We set $N_1 = G_1 \cap N$, $N_2 = G_2 \cap N$ and $N_{12} =
G_{12}\cap N$.

\begin{lemma} \label{Nstructure}Suppose that $G$ is a completion of the Goldschmidt amalgam $(G_1,G_2,G_{12})$ and that
 $N$ is a  normal subgroup of $G$. Then the following
hold.
\begin{enumerate}
\item $N_1$ and $N_2$ are not both non-trivial $2$-groups.

\item If $ O^2(G_i)\le N$ for some $i\in \{1,2\}$, then $G = NG_{3-i}$.

\item  If $\langle O^2(G_{1}),O^2(G_{2}) \rangle \le N$, then $G = G_{12}N$.

\item  If $\langle O^2(G_{1}),O^2(G_{2}) \rangle \le N$, then $N$ is a completion of the Goldschmidt amalgam
$(N_1,N_2,N_{12})$.
\end{enumerate}
\end{lemma}

\begin{proof} Suppose that both $N_1$ and $N_2$ are $2$-groups.
Then $$N_1 \le O_2(G_1)\cap N \le G_{12}\cap N \le N_2$$ and a
similar argument shows  $N_2 \le N_1$. Hence $N_1=N_2=1$ and (1)
is true.

Next suppose, without loss of generality, that $O^2(G_1) \le N$.
Then, as $G = \langle O^2(G_1),G_2\rangle$, we have that $G = G_2
N$ and  (2) holds. Similarly, as $G = \langle O^2(G_1),
O^2(G_2),G_{12}\rangle$, if $\langle O^2(G_1),O^2(G_2)\rangle \le
N$, then $G = G_{12}N$ and  (3) holds.

Continue to assume $N \ge \langle O^2(G_1),O^2(G_2)\rangle$. Then
$N_{12}\in \Syl_2(N_1)\cap \Syl_2(N_2)$ and so $[N_1:N_{12}]=
[N_2:N_{12}]=3$. Also $\langle N_1,N_2\rangle$ is normalized by
$G_{12}$ and so by considering orders we have $N =\langle
N_1,N_2\rangle$. Assume that $K$ is a normal subgroup of $\langle
N_1,N_2\rangle$ which is contained in $N_{12}$ and which is chosen
of maximal order. Let $x \in G_{12}$, then $K^x$ is normalized by
$N_1$ and $N_2$ and is contained in $N_{12}$. But then the maximal
choice of $K$ forces $K=K^x$ and consequently $K \normal G_{12}N
=G$. Hence $K=1$ and the proof of (3) is complete.
\end{proof}

Using more detailed structure of the Goldschmidt amalgams we can
say more about the situation in Lemma~\ref{Nstructure} (2).

\begin{lemma} \label{class1or2}Suppose that $G$ is a completion of the Goldschmidt amalgam $\mathcal A=(G_1,G_2,G_{12})$ and that
 $N$ is a  normal subgroup of $G$. If for some $i\in \{1,2\}$, $N_i \ge
O^2(G_i)$ and $N_{3-i}$ is a $2$-group, then $O^2(G_i) \cong 3$,
$\mathcal A$ is in Goldschmidt class $\GG_1$ or $\GG_2$ and
$N_{12}$ is elementary abelian.
\end{lemma}

\begin{proof} Suppose that $O^2(G_i) \not\cong 3$. Then
$O_2(O^2(G_i))\not=1$ and $$O_2(O^2(G_i)) \le N_{12}= N_{3-i}\le
O_2(G_{3-i}).$$  Lemma~\ref{gold1}(4) shows that $\mathcal A$ is
in class $\GG_1$ or $\GG_2$ and then Lemma~\ref{gold1} (1) and (2)
show that $O^2(G_i) \cong 3$, a contradiction. Therefore
$O^2(G_i)\cong 3$ and Lemma~\ref{gold1} (1) and (2) again shows
that $\mathcal A$ is in class $G_1$ or $G_2$. Finally, as
$N_{12}\le O_2(G_{3-i})$, we have $N_{12}$ is elementary abelian.
\end{proof}

The next  lemma  helps to decide whether or not $\G(G,\AA)$ is
symmetric or semisymmetric.

\begin{lemma}\label{notsymmetric} Suppose that $H$ is a group and  $G\le H$ is a completion
of an amalgam $\AA=(G_1,G_2,G_{12})$. If there exists $x \in H$
such that $G_1^x= G_2$ and $x^2 \in G_{12}$, then $\Gamma(G,\AA)$
is a symmetric graph.
\end{lemma}

\begin{proof} Define an automorphism of $\G =\G(\G,\AA)$ by
$G_ih \mapsto G_{3-i}xh$ for $i=1,2$. Since $G_2 = G_1^x$ and $x^2
\in G_{12}$, this map is a graph automorphism. It follows that
$\Aut(\G)$ acts transitively on the vertices of $\G$ and thus $\G$
is symmetric.
\end{proof}

We use Lemma~\ref{notsymmetric} and, especially, Tutte's Theorem
in its pushing-up disguise in the next lemma.

\begin{lemma} \label{SylowComp}Suppose that $\AA =(G_1,G_2,G_{12})$ is a Goldschmidt amalgam and $G$ is a completion of
$\AA$. If $G_1$ and $G_2$ are both maximal subgroups of $G$, then
either
\begin{enumerate}
\item $\G(G,\AA)$ is symmetric; or
\item $G$ is a Sylow completion of $\AA$.
\end{enumerate}
\end{lemma}

\begin{proof} Suppose that $G$ is not a Sylow completion of $\AA$
and let $T > G_{12}$ be a $2$-group with $[T:G_{12}]=2$. Then for
$x \in T\setminus G_{12}$ we have  $x^2 \in G_{12}$. If $G_1^x =
G_2$, then Lemma~\ref{notsymmetric} implies $\Gamma(G,\AA)$ is
symmetric. So assume $G_1^x \not= G_2$. Since $G_1$ and $G_2$ are
both maximal subgroups of $G$, we get  $\langle G_1,T\rangle =
\langle G_2,T\rangle = G$. Therefore $G$ is a completion of the
amalgams $(G_1,T,G_{12})$ and $(G_2,T,G_{12})$. It follows from
Tutte's Theorem \cite{Tutte} that $G_1$ is isomorphic to one of
$3$, $\Sym(3)$, $\Sym(4)$ or $\Sym(4)\times 2$ and that $G_2 \cong
G_1$ (by orders and symmetry). Suppose that $G_1 \cong \Sym(4)$ or
$\Sym(4) \times 2$. For $i=1,2$, set $Q_i = O_2(G_i)$ and set $Q_3
= Q_1^x$. Then $Q_1, Q_2$ and $Q_3$ are all elementary abelian of
order $2^2$ when $G_1\cong \Sym(4)$ and $2^3$ when $G_1 \cong
\Sym(4)\times 2$; furthermore, they are all contained in $G_{12}$.
However $G_{12}$ has only two such subgroups. Hence $Q_3 \in
\{Q_1,Q_2\}$. Then, as, for $i=1,2$, $G_i = N_G(Q_i)$ we get that
$G_1^x \in \{G_1,G_2\}$, which is a contradiction. Thus $G_1 \cong
3$ or $\Sym(3)$. In the first case, as $G$ has even order and
$G_1$ is a maximal subgroup of $G$,  $G_1 \in \Syl_3(G)$ and $G_1
= Z(N_G(G_1))$. Hence $G$ has a normal $3$-complement by
Burnside's Normal $p$-complement Theorem \cite[7.2.1]{StelKur}.
Hence $G = G_1R$ where $R$ is a normal $3^\prime$-subgroup of $G$.
Furthermore, as $G_1$ is a maximal subgroup of $G$, $C_R(G_1)= 1$
and so $R$ is nilpotent by Thompson's Theorem~\cite[Theorem
10.2.1]{Gorenstein}. But then the maximality of $G_1$  and $T>1$
forces $G \cong \Alt(4)$. Hence we calculate that $\Gamma(G,\AA)$
is symmetric. Next assume that $G_1 \cong \Sym(3)$. If $G$ is not
a simple group, then an argument like the one just explained
results in $\Gamma(G,\AA)$ being symmetric. Hence we may suppose
that $G$ is a simple group. The maximality of $G_1$ in $G$ means
that $O_3(G_1)$ is self centralizing in $G$, hence we may apply a
Theorem of Feit and Thompson \cite{FeitThompson2} to get that $G
\cong \Alt(5)$. Again we calculate that $\G(G,\AA)$ is symmetric.
This final contradiction shows that $G$ is a Sylow completion of
$\AA$.
\end{proof}

The above lemma can easily be adapted to  show that if $G$ is a
completion of a Goldschmidt amalgam $\AA=(G_1,G_2,G_{12})$ and,
say, $G_1$ is maximal in $G$ but $\Gamma(G,\AA)$ is  not
symmetric, then $G$ is either  a Sylow completion of $\AA $  or
$\AA$ is of type $\GG_2^1$, $\GG_2^2$, or $\GG_{2}^4$.  There are
numerous examples of this situation to be found in the groups
$\PSL_2(p)$ with $p$ odd (see \cite{Lip}). Indeed these groups and
groups closely related to them  probably provide the only examples
of this phenomenon.

\section{Some characterizations of simple groups}

In this section we gather together results which are either
characterisations of simple groups by the structure of their Sylow
$2$-subgroups  or statements asserting the non-simplicity of
groups with certain isomorphism types of Sylow $2$-subgroup.

\begin{theorem} \label{simpleSylows}Assume  that $G$ is a non-abelian finite simple
group and $S \in \Syl_2(G)$.
\begin{enumerate}
\item Suppose that $S$ is  abelian. Then either  \begin{enumerate} \item $G\cong\PSL_2(2^n)$, $n \ge
2$;\item $G\cong\PSL_2(q)$, $q\equiv 3,5 \pmod 8$, $q \ge 11$;
\item $G\cong \J_1$; or
\item $G\cong {}^2\GG_2(3^{2n+1})$, $n \ge 2$.
\end{enumerate}In particular, $S$ is elementary abelian.
\item Suppose that $S$ is a dihedral group of order $8$. Then
either
\begin{enumerate} \item
$G\cong\PSL_2(q)$, $q\equiv 7,9 \pmod {16}$; or  \item
$G\cong\Alt(7)$.
\end{enumerate}
\item Suppose that $S$ is isomorphic to a Sylow $2$-subgroup of
$\PSU_3(3)$. Then either
\begin{enumerate}
\item $G \cong \PSL_3(q)$, $q \equiv 5 \pmod 8$; or \item $G \cong
\PSU_3(q)$, $q \equiv 3 \pmod 8$.
\end{enumerate}
\item Suppose that $S$ is isomorphic to a Sylow $2$-subgroup of
$\M_{12}$. Then either
\begin{enumerate}
\item $G \cong \GG_2(q)$, $q \equiv 3,5 \pmod 8$;  \item $G \cong
{}^3\DD_4(q)$, $q \equiv 3,5 \pmod 8$; or \item $G \cong \M_{12}$.
\end{enumerate}
\end{enumerate}
\end{theorem}
\begin{proof} Part (1) is proved in \cite{Bender}. See also
\cite[Theorem 16.6]{Gorenstein}. For (2) we refer to
\cite{GW1,GW2,GW3,GW4} and to \cite[Theorem 16.3]{Gorenstein} for
a straight forward statement of the result. Part (3) is the
combined result of \cite{B1,B2}. Part (4) follows from \cite{GH}
and by consideration of the order of a Sylow $2$-subgroup of the
group. We also remark that the Sylow $2$-subgroups in parts (2),
(3) and (4) are generated by less than $5$ elements. Hence the
results also follows from \cite{GH2}, though of course this result
also relies on the results just cited above.
\end{proof}

\begin{lemma} \label{2D8}Suppose that $G$ is a finite group and $S \in
\syl_2(G)$. If $S \cong 2 \times \Dih(8)$, then $G$ is not simple.
\end{lemma}

\begin{proof} Suppose that $G$ is simple.
We first of all note that the automorphism group of $S$ is a
$2$-group. It follows that $SC_G(S)=N_G(S)$. In particular, by
Burnside's Theorem \cite[Theorem 7.1.1]{Gorenstein}, no two
distinct involutions in $Z(S)$ are conjugate in $G$. Let $z \in
(S^\prime)^\#$. Then by Glauberman's $Z^*$-Theorem \cite{Z*},
there is a conjugate $x$ of $z$ in $S\setminus \langle z \rangle$.
Since no two of the involutions in $Z(S)$ are conjugate, $x \not
\in Z(S)$. Hence there exists $S_0\le S$  such that $x\in S_0$ and
$S_0 \cong \Dih(8)$. Then $S_0$ contains three $S_0$-conjugacy
classes of involutions. Let $z_1$ and $z_2$ be the two involutions
from $Z(S)$ which are not contained in $Z(S_0)$. Then by
Thompson's Transfer Lemma \cite[12.1.1]{StelKur}, $z_1$ and $z_2$
are conjugate to elements of $S_0$. Since $z_1$ and $z_2$ are not
conjugate and neither of them is conjugate to $x$, we have a
contradiction.
\end{proof}

Suppose that $\AA=(G_1,G_2,G_{12})$ is an amalgam of type
$\GG_5^1$ and let $\AA^*=(G_1^*,G_2^*,G_{12}^*)$ be the subamalgam
of $\AA$ of type $\GG_5$ (see Lemma~\ref{subamalgam}). Then
$G_{12}$ is isomorphic to a Sylow $2$-subgroup of $\Aut(\M_{12})$
and $G_{12}^*$ is isomorphic to a Sylow $2$-subgroup of $\M_{12}$.
For $i=1,2$, set $Q_i = O_2(G_i)$, $Z_i = \Omega_1(Z(Q_i))$. Also
set $V_2 = \langle Z_1^{G_2}\rangle$, $U_1 = \langle
V_2^{G_1}\rangle$ and $W_1 = [U_1,O^2(G_1)]$. The following
results are obtained by direct calculations in the $\GG_5^1$
amalgam. (They were performed using {\sc Magma}\cite {Magma}.)

\begin{enumerate}
\item $Z_1 \cong 2^2$, $U_1$ has order $2^5$, $W_1 \cong 4\times 4$ and $W_1$ is a characteristic subgroup of $Q_1$.
\item $Z_2 \cong 2$ and $V_2$ is elementary abelian of order $2^3$. \item Let $F_1 \le U_1$ be the unique normal subgroup of $G_1$ of order $2^3$.
Then $F=C_{G_1}(F_1)$ is elementary abelian of order $2^4$ and is
the unique elementary abelian subgroup of order $2^4$ in $G_{12}$.
 \item If $x \in
F\setminus G_{12}^*$, then $C_{G_{12}}(x)=F$. \item Every
involution of $G_{12} \setminus G_{12}^*$ is in $F$.
 \item $[F,G_{12}]$ is a normal subgroup of
$G_1$ and has order $8$.\item If $y \in [F,G_{12}]$, then
$C_{G_{12}}(y)$ is non-abelian.
 \item $[F,G_{12},G_{12},G_{12}]= Z_2$. \item The
orbits of $G_{1}$ on $F$ have lengths $8$, $4$, $3$ and $1$. \item
The amalgam $\AA_1$ has at most two conjugacy classes of elements
of order $2$. Both classes are represented by elements in
$[F,G_{12}]$.
\end{enumerate}

\begin{lemma} \label{AutM12Sylows}Suppose that $G$ is a Sylow completion of an amalgam
of type $\GG_5^1$. Then $G$ is not  a simple group.
\end{lemma}

\begin{proof} Let $S = G_{12}$ and $x \in F \setminus G_{12}^*$.
Suppose that $x$ is conjugate to an element $y$ of $G_{12}^*$.
Then, from the fusion of involutions in $\AA_1$ as described in
(10), we may assume  that $y\in F$. Since $x$ is conjugate to $y$
and, by (7), $C_S(y)$ is non-abelian, if $R \in \Syl_2(C_G(x))$,
then $R$ is non-abelian. Select  $R\in \Syl_2(C_G(x))$ with $F \le
R$ and let $T \ge R$ be a Sylow $2$-subgroup of $G$ which contains
$R$. Then $F $ is normal in $T$ by (3).  Set $X = \langle
G_{1},T\rangle$. Then $X \le N_G(F)$. The orbits of $G_1$ on $F$
have lengths $8$, $4$, $3$ and $1$ by (9) and $T$ has orbits of
lengths $8$, $4$, $2$, $1$ and $1$ by (8). Since $C_T(x)
>F$, we have that $x$ not contained in the $T$-orbit of length
$8$. Therefore $T$ does not preserve  the orbit of $G_1$ of length
$8$. Since $11$ does not divide $|\GL_4(2)|$, we conclude that $X$
is either transitive on $F^\#$, or it has an orbit of length $3$
and an orbit of length $12$. Consider the latter case. Then $Z_1$
is normal in $X$. Let $Y$ be the subgroup of $X$ which centralizes
both $Z_1$ and $F/Z_1$ (so $X/Y $ is isomorphic to a subgroup of
$\Sym(3) \times \Sym(3)$). Since $O^2(G_1)$ centralizes $F/Z_1$
and $X$ has an orbit of length $3$ on $F/Z_1$, we have that $9$
divides $|X/Y|.$ It follows that $|C_X(Z_1)Y/Y|$ is also divisible
by $3$. Notice that $Q_1 \in \syl_2(Y)$. By the Frattini argument
we can choose $K \in \syl_3(C_X(Z_1))$ such that $K \le N_G(Q_1)$.
Now $K$ centralizes $Z_1$ and, as, by (1), $W_1$ is characteristic
in $Q_1$, $K$ normalizes $W_1$. Since, also by (1), $W_1 \cong
4\times 4$, we have that $K$ centralizes $W_1$. Note that $[F,K]$
has order $4$ with $[F,K]\cap Z_1=1$. On the other hand, $[F,K]$
is normalized by $Q_1=FW_1$ and so $[F,K]\cap Z_1>1$ by the
definition of $Z_1$.  It follows that $K $ centralizes $F$ and we
have a contradiction. Hence we may, and do, assume that $X$
operates transitively on $F^\#$. Therefore, upon  setting $Y=
C_G(F)$, $X/Y$ is isomorphic to a subgroup of $\GL_4(2)$ which
acts transitively on $F^\#$ and has Sylow $2$-subgroups isomorphic
to $S/F\cong \Dih(8)$. As $X$ acts transitively on $F^\#$, we also
have $O_2(X/Y)=1$. It follows that $X/Y$ is isomorphic to one of
the following groups: $\Sym(5)$, $\Alt(6)$, $\Alt(7)$, $(3 \times
\Alt(5)):2$. We now identify $\GL_4(2)$ with $\Alt(8)$  to
describe these subgroups. (Note that this isomorphism sends
elements of order $3$ with fixed points of $F^\#$ to elements of
cycle shape $3^2$.) Since $O^2(G_1)$ normalizes $Z_1$ and
centralizes $F/Z_1$, we see that $G_1Y/Y \cong \langle (123)
(456),(132)(567),(12)(45)\rangle$. Hence $X/Y$ is either
$\Alt(7)$, $(3 \times \Alt(5)):2$. Let $K$ be such that $K/Y$ has
order $3$ and $K/Y$ is normalized by $SY/Y$. Let $K_3 \in
\syl_3(K)$. Then $N_X(K_3) \cap F=1$ and $N_{KS}(K_3)Y = KS$.
Hence $N_X(K_3)$ has Sylow $2$-subgroups isomorphic to $\Dih(8)$.
It follows from Gasch\"utz's Theorem \cite[3.3.2]{StelKur} that
$S$ is isomorphic to a split extension $2^4:\Dih(8)$. But then $S$
has two elementary abelian subgroups of order $2^4$, in
contradiction to (3). This contradiction shows that $x$ cannot be
conjugate to an involution in $G_{12}^*$ and the  Thompson's
Transfer Lemma \cite[12.1.1]{StelKur} implies that $G$ has a
subgroup of index $2$. In particular, $G$ is not a simple group.
\end{proof}

\section{Sylow completions of Goldschmidt amalgams}

Throughout this section we assume that $\mathcal
A=(G_1,G_2,G_{12})$ is a Goldschmidt amalgam and that $G$ is a
Sylow completion of $\mathcal A$. We recall from the introduction
that a normal subgroup $R$ of $G$ is called \emph{regular} if it
acts semiregularly on $\G(G,\AA)$. Of course this is the same as
saying that $R_1 = R_2=1$. If $R$ is a regular normal subgroup of
$G$, then, as $G$ is a Sylow completion of $\mathcal A$, $R$ has
odd order and is consequently soluble by the Odd Order Theorem
\cite{FeitThompson}. The first lemma is trivial to prove.

\begin{lemma}\label{regular1} Suppose that $R$ is a regular normal
subgroup of $G$. Then $G/R$ is a Sylow completion of $\mathcal
A$.\qed
\end{lemma}

\begin{lemma} \label{regnormal}Suppose that $M$ and $N$ are regular normal
subgroups of $G$. Then either
\begin{enumerate}
\item $MN$ is a regular normal subgroup of $G$; or

\item $G$ is soluble.

\end{enumerate} In particular, if $G$ is not soluble, then $G$ the
product of all the regular normal subgroups of $G$ is the unique
maximal regular normal subgroup.
\end{lemma}

\begin{proof} Suppose that $M$ and $N$ are regular normal subgroups of $G$. Then, by Lemma~\ref{Nstructure} (1), $M$
and $N$ and consequently $MN$ have odd order. Hence $MN$ is
soluble. If $MN$ is not regular, then, by Lemma~\ref{Nstructure}
(1), either $O^2(G_1) \le MN$ or $O^2(G_2) \le MN$. It follows
that $G = G_1MN$ or $G = G_2 MN$. In both cases $G$ is soluble.
\end{proof}

\begin{lemma} \label{soluble} Assume that $G$ is soluble. Then
there exists a regular normal subgroup $R$ of $G$ such that one of
the following holds.

\begin{enumerate}
\item $G/R \cong 3$ or $\Sym(3)$ and $\mathcal A$ is of
type $\GG_1$ or $\GG_1^1$ respectively; or
\item$G/R \cong 3^2$, $3^2.2$,  $3 \wr 2$  or $\Sym(3) \times \Sym(3)$ and $\mathcal A$ is of type $\GG_1$,
$\GG_1^1$, $\GG_1^2$ or $\GG_1^3$ respectively; or
\item $G/R \cong 3^3.\Alt(4)$,
$3^3.\Sym(4)$, $3^3.\Sym(4)$, $ \Sym(3)\wr 3$ or $\Sym(3)\wr
\Sym(3)$ and $\mathcal A$ is of type $\GG_2$,  $\GG_2^1$,
$\GG_2^2$, $\GG_2^3$ or $\GG_2^4$ respectively.
\end{enumerate}
\end{lemma}

\begin{proof} Suppose that $N$ is a normal subgroup of $G$ of minimal order such
that either $N_1 \ge O^2(G_1)$ or $N_2 \ge O^2(G_2)$ (or both).
Assume without loss of generality that $O^2(G_1) \le N$.  Choose a
normal subgroup $R$ of $G$ contained in $N$ such that $N/R$ is a
minimal normal subgroup of $G/R$. Lemma~\ref{Nstructure} and the
minimal choice of $N$ shows that $R$ is a regular normal subgroup
of $G$. Also, as $O^2(G_1)\le N$ and $G$ is soluble, $N/R$ is an
elementary abelian $3$-group. In particular, $N_1$ is cyclic of
order $3$ and $G/N \cong G_2$. Since $G = \langle
O^2(G_1),G_2\rangle$ and $G_2= G_{12}T$ for $T \in \Syl_3(G_2)$,
we get that
$$N = \langle N_1R^{G_2}\rangle = \langle
N_1R^{G_{12}T}\rangle =\langle N_1^{T}R\rangle.$$ It follows that
$\langle N_1^{T}\rangle R/R = N/R$ has order dividing $27$.

Assume that $|N/R| = 3$. Then $O^2(G_1)$ and $O^2(G_2)$ commute
mod $R$. In particular,   both $O^2(G_1)$ and $O^2(G_2)$ are
cyclic of order $3$. Hence $\mathcal A$ is in class $\GG_1$ by
Lemma~\ref{gold1}. If $O^2(G_1)R=O^2(G_2)R$, then $|G/N| \le 2$,
$G_1 \cong G_2$ and the possibilities in (1) arise. If
$|O^2(G_1)O^2(G_2)R/R|= 3^2$, then the possibilities in (2) occur.

Suppose that $|N/R| \ge 3^2$. Let $M$ be the preimage of
$C_{G/R}(N/R)$. Then $M$ is normal in $G$ and, as $O^2(G_2)$ does
not normalize $O^2(G_1)R/R$, we have that $M_2$ is a $2$-group. It
follows that $M_1 = O^2(G_1)M_2$ and that $M_2$ is normal in $G$.
But then  $M_2=1$ and so $M=N$. Therefore, $G/N$ acts faithfully
and irreducibly on $N/R$. In particular, $\AA$ is in class
$\GG_2$. If $|N/R|= 27$, then $G_2 \cong G/N$ is isomorphic to a
subgroup of $\GL_2(3)$ and we have a contradiction to the
structure of $G_2$. Thus $|N/R|= 27$ and we obtain the
possibilities listed in  (3).
\end{proof}

\begin{remark}  The groups $G/R$ and the corresponding amalgams described  in Lemma~\ref{soluble} (3) each
determine graphs isomorphic to the  the Gray graph, it is the
smallest semisymmetric graph \cite{Bouwer}.
\end{remark}

From here on we assume that $G$ is not soluble and we let $R$ be
the unique maximal regular normal subgroup of $G$ (see
Lemma~\ref{regnormal}). Let $N$ be a normal subgroup of $G$ with
$N/R$ a minimal normal subgroup of $G/R$. If $N$ is soluble, then,
because of  the maximality of $R$, Lemma~\ref{Nstructure} implies
that $G = G_1N$ or $G = G_2N$ and we have a contradiction. It
follows that $N/R$ is a direct product of non-abelian simple
groups and $R$ is the maximal soluble normal subgroup of $G$. We
fix this notation for the remainder of the paper.

\begin{lemma} \label{simple} The quotient  $N/R$ is a non-abelian simple group and $G/R$ is isomorphic to a subgroup of
$\Aut(N/R)$ containing $\Inn(N/R)$.
\end{lemma}

\begin{proof} Assume
that $N/R$ is not simple. Then $G>N$ and, as $G_{12}\in
\syl_2(G)$,  $N_{12} = G_{12} \cap N \in \syl_2(N)$. Since  a
Sylow $2$-subgroup of a non-abelian simple group has order at
least $4$ \cite[7.2.2]{StelKur}, we have  $|N_{12}| \ge 2^4$.
Using that fact that $|N_{12}| \ge 2^4$ and the structure of the
Goldschmidt amalgams, we see that $G_{12} \not \le N$. Hence
$|G_{12}| \ge 2^5$. So $\mathcal A$ is of type $\GG_4^1$ or
$\GG_5^1$ and $(N_1,N_2,N_{12})$ is a Goldschmidt amalgam of type
$\GG_4$ or $\GG_5$ by Lemma~\ref{subamalgam}. But then, from the
structure of these two  amalgams, we have that $Z(N_{12})$ is
cyclic, and of course this is incompatible with $N/R$ being a
direct product of more than one simple group. This contradiction
shows that $N/R$ is a simple group.

Let $C= C_G(N/R)$. Then $C\ge R$ and $C \normal G$. If $C
> R$, then the maximality of $R$ means that, say, $O^2(G_1) \le C$ and $O^2(G_2) \le
N$. But then $[O^2(G_1),O^2(G_2)]\le R$ and $K =
O^2(G_1)O^2(G_2)R$ is normal in $\langle
O^2(G_1),O^2(G_2),G_{12}\rangle= G$, which is a contradiction as
$G$ is not soluble. Hence $C=R$ and we have that $G/R$ is
isomorphic to a subgroup of $\Aut(N/R)$ containing $\Inn(N/R)$ as
claimed.
\end{proof}

We now examine the possibilities for $N/R$ by examining the
possibilities for the Sylow $2$-structure of $N$. In the next two
lemmas we frequently use well-known structural facts about
$\PSL_2(q)$, $q=p^a$ for a prime $p$. We refer the reader to the
theorem of Dickson which can be found in \cite[Hauptsatz
8.27]{Huppert}. In particular, we note that for odd primes $p$,
the centralizer of an involution in $X=\PSL_2(p^a)$ is dihedral of
order $p^a+1$ if $p^a\equiv 3 \pmod 4$ and order $p^a-1$ if
$p^a\equiv 1 \pmod 4$. We also note that if $F$ is a fours group
of $X$, then $N_X(F) \cong \Sym(4)$ if $q\not \equiv 3, 5 \pmod 8$
and $N_X(F) \cong \Alt(4)$ if $q\equiv 3, 5 \pmod 8$.

\begin{lemma} \label{abelian} Suppose that $N$ has abelian Sylow $2$-subgroups.
Then $(N_1,N_2,N_{12})$ is a Goldschmidt amalgam of type $\GG_1^3$
or $\GG_2$ and $ N/R\cong \PSL_2(p)$ with $p$ a prime, $p\equiv
11, 13 \pmod {24}$. Furthermore, either $G=N$ or $G/R \cong
\PGL_2(p)$ with $p$ a prime, $p\equiv 11,13 \pmod {24}$ and
$\mathcal A$ is of  Goldschmidt type $\GG_2^1$.

\end{lemma}

\begin{proof} We may  assume that $R=1$.   By Theorem~\ref{simpleSylows} (1), we have that $N
\cong \PSL_2(q)$ with $q\equiv 3,5\pmod 8$ and $q \ge 11$,
$\PSL_2(2^n)$ with $n \ge 2$,  ${}^2\GG_2(3^{2t+1})$ with  $t \ge
2$ or $\J_1$. Without loss of generality we assume that $N_{1} \ge
O^2(G_1)$.

Assume that $N_2 \in \Syl_2(N)$. Then, as $N_2 \le O_2(G_2)$, we
have that $\mathcal A$ is in class $\GG_1$ or $\GG_2$ by
Lemma~\ref{class1or2}. Furthermore, we have that $N_1 $ is
contained in the centralizer of an involution of $N$ and $4\le
|N_2|\le 8$. Assume that $|N_{2}| = 8$. Then $N_1 \cong \Sym(3)
\times 2 \times 2$. However, the centralizer of an involution in
$\J_1$  is isomorphic to $2 \times \Alt(5)$, in
${}^2\GG_2(3^{2t+1})$ is isomorphic to $2\times \PSL_2(3^{2t+1})$
and in $\PSL_2(2^3)$ is abelian. Since none of these groups
contains a subgroup isomorphic to $\Sym(3) \times 2 \times 2$, we
conclude that $|N_{12}| = 4$. It follows that $N \cong \PSL_2(q)$
with $q\equiv 3,5\pmod 8$ and $N_1 \cong 2 \times \Sym(3)$. Now,
from the structure of $\PSL_2(q)$, $X=N_{N}(N_{12}) \cong
\Alt(4)$. It follows that $N_{G}(N_{12}) = G_2 X$ and, since,
under our present hypothesis, $O^2(G_2) \not \le N$,  we have that
$9$ divides $|N_{G}(N_{12})|$. Therefore there is a non-trivial
element $x\in O_3(G_2 X)$ such that $[x,N_{12}]=1$. In particular,
$x$ normalizes $C_{N}(y)$ for  $y \in O_2(N_1 )^\#$ and as
$C_{N}(y)$ is dihedral and $x$ has order $3$, we infer that $x$
centralizes $N_1$. Therefore, $x$ centralizes $\langle
O^2(G_2),O^2(G_1)\rangle \ge N$ and this contradicts
Lemma~\ref{simple}. This contradiction shows that $N_2 \ge
O^2(G_2)$. Hence $\mathcal A_N= (N_1,N_2,N_{12})$ is a Goldschmidt
amalgam. From Table~\ref{GoldList} we read that $\mathcal A_N$ is
of type $\GG_1^3$, $\GG_2$ or  $\GG_2^3$. The latter possibility
is ruled out just as above as our candidates for $N$ do not
contain a subgroup isomorphic to $\Sym(3) \times 2 \times 2$.
Therefore, $N \cong \PGL_2(q)$ with $q\equiv 3,5\pmod 8$. We
choose notation so that $N_1 \cong 2\times \Sym(3)$. Let $p$ be a
prime chosen so that $q= p^n$ for some $n$. Then as $3$ divides
the order of the centralizer of an involution in $N$ we have that
$p \not =3$. Furthermore, we have that $p\equiv 3, 5 \pmod 8$ and,
by Sylow's Theorem we may assume that $N_{12} \le H \le N$ where
$H \cong \PSL_2(p)$. Recalling that the centralizer of an
involution in $\PSL_2(q)$ is a dihedral group of order $q-1$ if $q
\equiv 5 \pmod 8$ or order $q+1$ if $q\equiv 3 \pmod 8$, we
require $q \equiv 11, 13 \pmod {24}$. It follows that $p\equiv 11,
13 \pmod {24}$. Hence, as dihedral groups of order divisible by
$3$ contain a unique cyclic subgroup of order $3$, we have $N_1
\le H$. Now $N_H(N_{12})= N_N(N_{12})\cong \Alt(4)$ and so it
follows that $\langle N_1,N_2\rangle = H$. Thus $N = H$ and we
have the first of our claims. Now assume that $G > N$. Then by
considering Table~\ref{GoldList} we deduce the stated structure of
$G$ and type of $\mathcal A$.
\end{proof}

We remark that one of the consequences of  Lemmas~\ref{class1or2}
and \ref{abelian}  is that $(N_1,N_2,N_{12})$ is a Goldschmidt
amalgam whenever $G$ is not soluble. We call this subamalgam
$\mathcal A_N$.

\begin{lemma}\label{dihedral} Suppose that $N_{12}\cong \Dih(8)$.
Then
\begin{enumerate}
\item $\mathcal A_N$   is of type $\GG_2^2$, $N/R \cong \Alt(7)$ and either  $G=N$ or  $\mathcal A $ is
of type $\GG_2^4$ and $G/R \cong \Sym(7)$.

\item $\mathcal A_N$ is of type $G_3$ and $N/R \cong \PSL_2(p)$ or $\PSL_2(p^2)$
with $p$ a prime satisfying  $p \equiv 7,9 \pmod {16}$ in the
former case and $p\equiv 3,5 \pmod 8$  in the latter case. If $G
>N$, then $\mathcal A$ is of type $\GG_3^1$ and $G/R \cong
\PSigmaL_2(p^2)$ with $p$ a prime satisfying $p \equiv 3,5 \pmod
8$.

\item $\mathcal A_N$ is of type $G_2^1$ and $N/R \cong \PSL_2(p)$ or $\PSL_2(p^2)$ with $p$ a prime
satisfying $p \equiv 23, 25 \pmod {48}$ in the former case and
$p\equiv 5,19 \pmod {24}$  in the latter case.  If $G >N$, then
$\mathcal A$ is of type $\GG_2^4$ and $G/R \cong \PSigmaL_2(p^2)$
with $p$ a prime satisfying $p \equiv 5,19 \pmod {24}$.
\end{enumerate}
\end{lemma}

\begin{proof} We assume that $R=1$. From Table~\ref{GoldList}, we see $\mathcal A_N$ is of type $\GG_2^1$, $\GG_2^2$ or
$\GG_3$. Since $N$ is simple, Theorem~\ref{simpleSylows} (2)
implies that $N \cong \PSL_2(q)$, $q \equiv 7,9 \pmod {16}$ or $N
\cong \Alt(7)$. In the latter case it is straightforward to show
that (1) holds.

So suppose that $q = p^a$ with $p$ a prime and that $N \cong
\PSL_2(q)$, $q \equiv 7,9 \pmod {16}$. If $p \equiv 3,5 \pmod 8$,
then $q$ is a square and we let $H$ be a subgroup  of $N$
isomorphic to $\PSL_2(p^2)$ and if $p \equiv   1,7 \pmod 8$, we
let $H$ be a subgroup  of $N$ isomorphic to $\PSL_2(p)$. Then $H$
has dihedral Sylow $2$-subgroups of order $8$. By Sylow's Theorem,
we may suppose that $H$ contains $N_{12}$. Let $F_1$ and $F_2$ be
the distinct elementary abelian subgroups of order $4$ in $N_{12}$
and, for $i=1,2$, put $X_i= N_G(F_i)$. Then $X_1\cong X_2 \cong
\Sym(4)$ and $H = \langle X_1, X_2\rangle$. In particular, if
$\AA_N$ is of type $\GG_3$, then $N_i=X_i$ and we have that $N =
H$. Let us continue to assume that $\AA_N$ has type $\GG_3$ and
additionally suppose that $G > N$. Then $\AA$ has type $\GG_3^1$.
If $p \equiv 1,7 \pmod 8$, then $G$ is isomorphic to $\PGL_2(p)$
and so has $G$  dihedral Sylow $2$-subgroups, a contradiction as
$G_{12} \cong 2 \times \Dih(8)$. Therefore, $p\equiv 3,5 \pmod 8$
and $H \cong \PSL_2(p^2)$. Then $G$ is isomorphic to a subgroup of
$\PGammaL_2(p^2)$. It is now easy to see that $G \cong
\PSigmaL_2(p^2)$.

We now consider the possibility that $\mathcal A_N$ is of type
$\GG_2^1$ or $\GG_2^2$. The latter case is impossible as the
centralizer of an involution in $\PSL_2(q)$ is dihedral and so
does not contain a subgroup isomorphic to $2^2\times 3$. So assume
$\mathcal A_N$ is of type $\GG_2^1$. We have that $N \cong
\PSL_2(q)$ with $q=p^n \equiv 7, 9\pmod {16}$. We additionally
require that the centralizer of an involution contains a subgroup
isomorphic to $\Dih(24)$. Therefore $q \equiv 23, 25 \pmod {48}$.
If $p\equiv 23, 25 \pmod {48}$, then just as above we see that $N
\cong \PSL_2(p)$. So assume that $p\equiv 3,5 \pmod 8$. Then $q$
is a square. If $p \equiv 11, 13 \pmod {24}$, then there is a
subgroup $K$ of $N$ with $K \cong \PGL_2(p)$ such that $K= \langle
N_1,N_2\rangle$ and we have a contradiction as $K <N$. Therefore,
$p \equiv 5, 19 \pmod {24}$ and the configurations in (3) arise.
\end{proof}

\begin{lemma} \label{GG4} Suppose that $\AA_N$ is of type $\GG_4$. Then $N/R \cong
\PSL_3(p)$ with $p$ a prime such that $p \equiv 5 \pmod 8$ or
$N/R\cong \PSU_3(p)$ with $p$ a prime such that $p \equiv 3 \pmod
8$. Furthermore, if $G >N$, then $\AA$ is of type $\GG_4^1$ and
$G/R \cong \PSL_3(p).2$ when $p\equiv 5 \pmod 8$ or $\PSU_3(p).2$
when $p\equiv 3\pmod 8$.
\end{lemma}

\begin{proof} Since in amalgams of type $\GG_4$ we know that $G_{12}$ is isomorphic to a Sylow $2$-subgroup
of $\PSU_3(3)$,  Theorem~\ref{simpleSylows} (3) implies that $N/R
\cong \PSU_3(q)$ with $q \equiv 3 \pmod 8$ or $N/R \cong
\PSL_3(q)$ with $q\equiv 5 \pmod 8$. Now applying \cite[Theorem
1.1]{ParkerRowley3} we have that $q$ is prime. This concludes the
proof of the lemma.
\end{proof}

\begin{lemma} \label{G5structure} Suppose that $(G_1,G_2,G_{12})$
is a Goldschmidt amalgam of type $\GG_4^1$ or $\GG_5$ (with the
orientation of $G_1$ and $G_2$ as in Table~\ref{GoldList}). Then
\begin{enumerate}
\item $|Z(G_{12})|=2$;
\item $\Omega_1(Z(O_2(G_1)))= \Omega_1([G_{12},G_{12}])$.
\item $O_2(G_2)$ is the unique extraspecial subgroup of $G_{12}$ of
order $2^5$.
\end{enumerate}
\end{lemma}

\begin{proof} Since in both these amalgams
$G_{12}$ is isomorphic to a Sylow $2$-subgroup of $\M_{12}$, this
result  is easily calculated using {\sc Magma} \cite{Magma} for
example.\end{proof}

\begin{lemma}\label{GG5} Suppose that $\mathcal A_N$ is of type $\GG_4^1$ or $\GG_5$. Then
 $\mathcal A_N$ is of type $\GG_5$ and one of the following
holds.
\begin{enumerate}
\item $N/R \cong \M_{12}$; or
\item $N/R \cong \GG_2(p)$ with $p$ a prime such that $p \equiv 3,5 \pmod
8$.
\end{enumerate}
Furthermore, if $G >N$, then $G$ is of type $\GG_5^1$ and $G/R
\cong \Aut(\M_{12})$ or $\Aut(\GG_2(3))$.
\end{lemma}

\begin{proof} As always  we assume that $R=1$. In both amalgams we have that $N_{12}$ is
isomorphic to a Sylow $2$-subgroup of $\M_{12}$. Hence
Theorem~\ref{simpleSylows} implies that $G$ is isomorphic to
$\M_{12}$, $\GG_2(q)$ or ${}^3\DD_4(q)$ with $q = p^a$ and
$q\equiv p \equiv 3,5 \pmod 8$. If $N\cong \M_{12}$, we  use
\cite{Atlas} to verify that $N$ is a Sylow completion of $\mathcal
A_N$. So assume that $N\cong \GG_2(q)$ or ${}^3\DD_4(q)$ with $q
\equiv 3,5 \pmod 8$. By Kleidman \cite{Kliedman1,Kliedman2} $N$
contains a subgroup $H$ isomorphic  to $\GG_2(p)$. By considering
orders, we have that $H$ contains a Sylow $2$-subgroup of $N$ and
so we may assume that $N_{12}\le H$. Let $Z = Z(N_{12})$. Then
$|Z|=2$ and $C_H(Z) \cong (\SL_2(p)\circ \SL_2(p)):2$ is contained
in $C_N(Z)$ where $C_N(Z)\cong (\SL_2(q)\circ \SL_2(q)):2$ when $N
\cong \GG_2(q)$ and $C_N(Z)\cong (\SL_2(q^3)\circ \SL_2(q)):2$
when $N \cong {}^3\DD_4(q)$. Since, by Lemma~\ref{G5structure},
$O_2(N_2)$ is uniquely determined in $N_{12}$, we can identify
$O_2(N_{12})\cong 2^{1+4}_+$ as a Sylow $2$-subgroup of
$O^2(C_H(Z))\cong \SL_2(p)\circ\SL_2(p)$. Then $$N_G(O_2(N_2)) =
N_{O^2(C_H(Z))}(O_2(N_2))N_{12} \sim 2^{1+4}_+.3^2.2.$$ In
particular, we see that $N_{G}(O_2(N_2))\le H$.  In $C_H(Z)$, the
two subgroups isomorphic to $\SL_2(p)$ are in fact normal
subgroups of $C_H(Z)$ (because they are generated by root groups
corresponding to long and short roots) and so
$O_3(N_G(O_2(N_2))/O_2(N_2))$ is inverted by $N_{12}/O_2(N_{12})$
(see \cite[Table 4.5.1]{GLS3} for a precise description of
$C_H(Z)$). It follows that there are exactly four subgroups $X_1$,
$X_2$, $X_3$ and $X_4$ of $N_G(O_2(N_2))$ each containing $N_{12}$
at index $3$ and each having $X_i/O_2(X_i) \cong \Sym(3)$ for
$i\in \{1,2,3,4\}$. Furthermore, we may assume that
$\e(X_1,O_2(X_1))= \e(X_2,O_2(X_2))=1$ and
$\e(X_3,O_2(X_2))=\e(X_4,O_2(X_2))=2$. Therefore, if $\mathcal
A_N$ is of type $G_4^1$, then $N_2\in \{X_1,X_2\}$ and, if
$\mathcal A_N$ is of type $G_5$, then $N_2\in \{X_3,X_4\}$.

By Lemma~\ref{G5structure}, $Z_1=\Omega_1([N_{12},N_{12}])$ and so
is uniquely determined by $N_{12}$. Additionally $O_2(N_1) =
C_{N_1}(Z_1) = C_{N_{12}}(O_2(N_1))$ is also uniquely determined.
Therefore $N_1 \le N_G(O_2(N_1))\le  N_G(Z_1)$. By
\cite[Proposition 2.4]{Kliedman1} and \cite[Proposition 2.6 (case
4)]{Kliedman2}
$$N_G(Z_1) \sim \begin{cases} (\mathbb
Z_{(q+1)}\times \mathbb Z_{(q+1)}).\Dih(12)& N \cong \GG_2(q), q
\equiv 3 \pmod 8\cr (\mathbb Z_{(q-1)}\times \mathbb
Z_{(q-1)}).\Dih(12)& N \cong \GG_2(q), q \equiv 5 \pmod 8\cr
(\mathbb Z_{(q^3+1)}\times \mathbb Z_{(q+1)}).\Dih(12)& N \cong
{}^3\DD_4(q), q \equiv 3 \pmod 8\cr (\mathbb Z_{(q^3-1)}\times
\mathbb Z_{(q-1)}).\Dih(12)& N \cong {}^3\DD_4(q), q \equiv 5
\pmod 8.\end{cases}$$  Let $K = C_{G}(Z_1)$. Then $N_N(Z_1)/K
\cong \Sym(3)$. Notice that $K$ has an abelian subgroup $T$ at
index $2$ with
$$T \cong
\begin{cases} \mathbb Z_{(q+1)}\times \mathbb Z_{(q+1)}& N \cong
\GG_2(q), q \equiv 3 \pmod 8\cr \mathbb Z_{(q-1)}\times \mathbb
Z_{(q-1)}& N \cong \GG_2(q), q \equiv 5 \pmod 8\cr \mathbb
Z_{(q^3+1)}\times \mathbb Z_{(q+1)}& N\cong {}^3\DD_4(q), q \equiv
3 \pmod 8\cr \mathbb Z_{(q^3-1)}\times \mathbb Z_{(q-1)}& N \cong
{}^3\DD_4(q), q \equiv 5 \pmod 8.\end{cases}.$$ We have
$N_N(Z_1)/T \cong \Dih(12)$. Now $O_2(N_1)$ induces a group of
order $2$ on $T$ and this action inverts every element of $T$. It
follows that $N_K(O_2(N_1)) = O_2(N_1)\cap T$. By the Frattini
argument, we have $N_{N_H(Z_1)}(O_2(N_1))K = N_H(Z_1)$ and so we
get $N_{N_H(Z_1)}(O_2(N_1))$ is a group of order $2^6.3$.
Therefore, $N_1 = N_{N_H(Z_1)}(O_2(N_1))$ and thus $N_1$ is
uniquely determined by $N_{12}$. We have shown that there are
exactly four possibilities for the amalgam $(N_1,N_2,N_{12})$,
namely $(N_1,X_i,N_{12})$ with $i\in \{1,2,3,4\}$ and all these
amalgams are contained in $H$. It follows that $N=H$. We refer
again to \cite[Theorem A]{Kliedman2} to get that the subgroups of
$H$ which contain $N_{12}$ are  $M_1\cong2^3\udot \GL_3(2) $
(non-split), $M_2\cong\GG_2(2) \cong \GG_2(\mathbb Z)$,
$M_3=N_H(Z_1)$, $M_4=C_{H}(Z)$, $M_5^+\cong\SL_3(p):2$ when $p
\equiv 5 \pmod 8$ and $M_5^-\cong\SU_3(p):2$ when $p \equiv 3
\pmod 8$. The group $M_1\cong 2^3\udot\GL_3(2)$ contains both
$N_1$ and one of  $X_3$ and $X_4$. Choose notation so that $X_3
\le M_1$. Then $M_1=\langle N_1,X_3\rangle$.  The group
$M_2\cong\GG_2(2)$ contains one of $X_1$ and $X_2$, so we may
suppose that $M_2= \langle N_1,X_1\rangle$. By definition $\langle
N_1,X_i\rangle$ is not contained in $M_3$ or $M_4$. When $p \equiv
3 \pmod 8$, $\langle N_1,X_2\rangle = \SU_3(p):2 \cong M_5^-$ and
when $\equiv 5 \pmod 8$, $\langle N_1,X_2\rangle \cong
\SL_3(p):2\cong M_5^+$. It remains to consider $\langle
N_1,X_4\rangle$. The groups $M_1$ and $M_2$ both contain only two
proper subgroups containing $N_{12}$ properly. So $X_4$ is not
contained in either of these subgroups. If $X_4 \le \M_5^\pm$,
then $\langle X_2,X_4\rangle \sim 2^{1+4}_+.3^2.2$ is also a
subgroup. But in $\SL_3(p)$, the centralizer of an involution is
isomorphic $\GL_2(p)$ and in $\SU_3(p)$ is isomorphic to
$\GU_2(p)$. We conclude that $X_4\not\le M_5^\pm$. The only
possibility that remains is that $\langle N_1,X_4\rangle = N$.
Hence $N$ is a completion of a Goldschmidt amalgam of type $G_5$
and we conclude that (2) holds. Finally, if $G >N$, then, as
$\Aut(\GG_2(p))
>\Inn(\GG_2(p))$  only for $p=3$ (see \cite{Atlas}) we infer that
if $G>N$, then $G\cong \Aut(\GG_2(3))$. If $N \cong \M_{12}$ we
use \cite{Atlas} to verify the result.
\end{proof}

We remark that the subgroups $2^3\udot\GL_3(2)$ and $\GG_2(2)$ are
critical in describing the GAB (geometry which is almost a
building) which has diagram $\begin{picture}(65,10)
\multiput(2,5)(30,0){3}{\circle*{4}} \put(2,5){\line(1,0){30}}
\put(32,6.8){\line(-1,0){30}} \put(32,3.2){\line(-1,0){30}}
\put(60,5){\line(-1,0){30}}
\end{picture}$ and was discovered by Kantor \cite{Kantor}.

We  now prove Theorem~\ref{MainTheorem}.

\begin{proof}[The Proof of Theorem~\ref{MainTheorem}]
Suppose that $G$ is a Sylow completion of the Goldschmidt amalgam
$(G_1,G_2,G_{12})$. If $G$ is soluble, then by Lemma~\ref{soluble}
we have that $G$ is as described in the first three divisions of
Table~\ref{Answers}. So we have that $G$ is not soluble. By
Lemma~\ref{regnormal} there is a unique maximal regular normal
subgroup $R$ of $G$ and by Lemma~\ref{simple} $G/R$ has a unique
minimal normal subgroup $N/R$ which is a non-abelian simple group.
If $N/R$ has abelian Sylow $2$-subgroups, then the structure of
$G/R$ is given in Lemma~\ref{abelian} and listed in segments four
and five of Table~\ref{Answers}. So we may assume that $N_{12}$ is
not abelian. Hence, by Lemmas~\ref{class1or2} and
\ref{Nstructure}, we have that $N/R$ is a completion of a
Goldschmidt amalgam. From the structure of the amalgams in
Table~\ref{GoldList} and because of Lemmas~\ref{2D8} and
\ref{AutM12Sylows}, we have that the Sylow $2$-subgroups of $N$
are isomorphic to either a dihedral group of order $8$, or $N$ is
a completion of an amalgam of type $\GG_4$, $\GG_4^1$ or $\GG_5$.
All these possibilities  are discussed in Lemmas~\ref{dihedral},
\ref{GG4} and \ref{GG5} and result is the groups described in the
final divisions of Table~\ref{Answers}.\end{proof}

Of course Theorem~\ref{MainTheorem2} is just another way  of
expressing Theorem~\ref{MainTheorem}.

\begin{proof}[The Proof of  Corollary~\ref{gammaR}] We let $A =
\Aut(\Gamma/R)$. For the groups in division 2 of
Table~\ref{Answers}, $\Gamma/R$ has six vertices and is a
bipartite, cubic graph and so we infer that $\Gamma/R \cong
\K_{3,3}$ and that $A\cong \Sym(3)\wr 2$. So $\Gamma/R$ is
symmetric. For division 3, we have that $\G/R$ has $54$ vertices
and is the so-called Gray graph and we have $A \cong \Sym(3)\wr
\Sym(3)$.

For the remaining cases we will be a little bit more
sophisticated. We  recall the fundamental Theorem of Tutte
\cite{Tutte} which says that if $\Theta$ is a symmetric, cubic
graph and $v$ is a vertex of $\Theta$, then the stabiliser of $v$
in the automorphism group of $\Theta$ is isomorphic to one of $3$,
$\Sym(3)$, $\Sym(3)\times 2$, $\Sym(4)$ or $\Sym(4)\times 2$.
Therefore for $\Gamma/R$ to be symmetric, we require that $G_1$
and $G_2$ both isomorphic to subgroups of one of these groups. It
follows immediately that for the graphs in divisions 5, 6, 7, 8,
11, 12, 13 and 14, we have that $A \cong \PGL_2(p)$, $\Sym(7)$,
$\PSL_2(p)$, $\PSigmaL_2(p^2)$, $\PSL_3(p).2$, $\PSU_3(p).2$,
$\Aut(\M_{12})$ and $\GG_2(p)$, $p\not=3$ and $\Aut(\GG_2(3))$
respectively and that in each case $\Gamma/R$ is semisymmetric. In
all the other cases  the graphs are symmetric. The automorphism
groups are, in divisions 4, 9 and 10, respectively, $\PGL_2(p)$,
$\PGL_2(p)$, and $\PGammaL_2(p^2)$ respectively.
\end{proof}

\begin{proof}[The Proof of  Corollaries~\ref{biprimitive} and
\ref{uniprimitive}] Note that by Lemma~\ref{SylowComp} we may
assume that $G$ is a Sylow completion of $\AA$. We work our way
through the divisions of Table~\ref{Answers}. Suppose first that
$G$ is in division $1$. Then $G= G_1 R = G_2R$ and we see that for
$G_1$ to be a maximal subgroup of $G$ we must have that $R$ is a
minimal normal subgroup of $G$. Then either $R $ is cyclic of
prime order $p$ or is elementary abelian of order $p^2$. Since in
this case $\Aut(G)$ possesses an involution which centralizes
$G_1$ and inverts $R$, we see that the graphs obtained from these
groups are all symmetric by using Lemma~\ref{notsymmetric}.

In all the remaining divisions of Table~\ref{Answers}  we have
that $G_1R < G$ and $G_2 R<G$. Thus for either of $G_1$ or $G_2$
to be a maximal subgroup of $G$ we must have that $R=1$. Therefore
the groups and graphs that we have to investigate are those listed
in Table~\ref{SS} which are not symmetric. For $\AA$ a Goldschmidt
amalgam we use the ordering of $G_1$ and $G_2$ as given in
Table~\ref{GoldList}.

For the group in division 3 of Table~\ref{SS} we have that $G_1$
is a maximal subgroup of $G$ and $G_2$ is not. This then is the
graph in Corollary~\ref{uniprimitive} (1). The groups in division
5 of Table~\ref{SS} have $G_1 \cong \Sym(4)$ and $G_2 \cong
\Dih(12)$. $G_1$ is a maximal subgroup of $G $ and $G_2$ is a
maximal subgroup if and only if $p\in \{11,13\}$. Thus
$\PGL_2(11)$ and $\PGL_2(13)$ appear in
Corollary~\ref{biprimitive} and the remaining groups are listed in
Corollary~\ref{uniprimitive} (2).

The group in division 6 of Table~\ref{SS} does not act primitively
on either orbit of $G$ on the vertices of $\Gamma$. Suppose that
$G$ is in division 7 of Table~\ref{SS}. Then again $G_1\cong
\Sym(4)$ is maximal in $G$ and $G_2 \cong \Dih(24)$ is only
maximal in $G$ if $p=23$. Thus $\PGL_2(23)$ is listed  in
Corollary~\ref{biprimitive} and the others are listed in
Corollary~\ref{uniprimitive} (3).

Consider next  division 8 of Table~\ref{SS}. So $G \cong
\PSigmaL_2(p^2)$, $G_1 \cong \Sym(4) \times 2$ and $G_2 \cong
\Sym(3) \times \Dih(8)$. Now in $\PSigmaL_2(25)$, $G_1$ is
contained in a subgroup isomorphic to $\Sym(5)\times 2 = \PGL_2(5)
\times 2$, thus $G_1$ is not maximal, whereas $G_2$ is a maximal
subgroup. Hence $\PSigmaL_2(25)$ appears in
Corollary~\ref{uniprimitive} (4) (and not in
Corollary~\ref{biprimitive}). For the remaining groups $G_1$ is
maximal in $G$ and $G_2$ is not, so they are in
Corollary~\ref{uniprimitive} (5).

Moving on to divisions 11 and 12 of Table~\ref{SS}, we quickly see
that the only possibilities are that $p=5$ in the former case and
that $p=3$ in the latter case. Thus $\PSL_3(5).2$ appears in
Corollary~\ref{uniprimitive} (6) and $\PSU_3(3).2 \cong \GG_2(2)$
is in Corollary~\ref{biprimitive}. Finally, $\Aut(\M_{12})$ is
primitive on both orbits so is in Corollary~\ref{biprimitive} and
in $G=\GG_2(p)$, as we saw in the proof of Lemma~\ref{GG5},
neither of $G_1$ nor $G_2$ is maximal in $G$. This completes the
proofs of Corollaries~\ref{biprimitive} and \ref{uniprimitive}.
\end{proof}

\end{document}